%% file: agt-4-33.tex
\documentclass{gtart_h}

\input agtout

\lognumber{33}
\volumenumber{4}
\volumeyear{2004}
\papernumber{33}
\pagenumbers{757}{780}
\received{11 December 2003} 
\accepted{23 August 2004}
\published{11 September 2004}

\usepackage{amssymb, amsmath,amscd} 
\usepackage[matrix,arrow]{xy}

\newtheorem{thm}{Theorem}
\newtheorem{lem}[thm]{Lemma}
\newtheorem{prop}[thm]{Proposition}
\newtheorem{cor}[thm]{Corollary}

\theoremstyle{definition}
\newtheorem{rem}[thm]{Remark}
\newtheorem*{ack}{Acknowledgements}

\renewcommand{\epsilon}{\ensuremath{\varepsilon}}
\renewcommand{\to}{\ensuremath{\longrightarrow}}
\newcommand{\rp}[1][2]{\ensuremath{{\mathbb R}P^{#1}}}
\newcommand{\Et}{\ensuremath{\mathbb E}^{\,2}}
\newcommand{\N}{\ensuremath{\mathbb N}}
\newcommand{\Z}{\ensuremath{\mathbb Z}}
\newcommand{\dt}{\ensuremath{\mathbb D}^{\,2}}
\newcommand{\St}[1][2]{\ensuremath{\mathbb S}^{#1}}
\newcommand{\FF}{\ensuremath{\mathbb F}}
\newcommand{\F}[1][2]{\ensuremath{\FF_{{#1}}}}
\newcommand{\sn}[1][n]{\ensuremath{S_{{#1}}}}
\renewcommand{\ker}[1]{\ensuremath{\operatorname{\text{Ker}}({#1})}}

\newcommand{\pnm}[1][n]{\ensuremath{P_{{#1}}(M)}}

\newcommand{\map}[4][\to]{\ensuremath{{#2} \co {#3} #1 {#4}}}
\newcommand{\id}{\ensuremath{\operatorname{\text{Id}}}}
\newcommand{\aut}[1]{\ensuremath{\operatorname{\text{Aut}}(#1)}}

\newcommand{\lhra}{\mathrel{\lhook\joinrel\to}}

\DeclareRobustCommand*{\up}[1]{\textsuperscript{#1}}
\renewcommand{\th}{\ensuremath{\up{th}}}
\newcommand{\ft}[1][n]{\ensuremath{\Delta_{#1}}}

\newcommand{\brak}[1]{\ensuremath{\{ #1 \}}}
\newcommand{\ang}[1]{\ensuremath{\left\langle #1\right\rangle}}

\newcommand{\reth}[1]{Theorem~\protect\ref{th:#1}}
\newcommand{\relem}[1]{Lemma~\protect\ref{lem:#1}}
\newcommand{\repr}[1]{Proposition~\protect\ref{prop:#1}}
\newcommand{\reco}[1]{Corollary~\protect\ref{cor:#1}}
\newcommand{\resec}[1]{Section~\protect\ref{sec:#1}}
\newcommand{\req}[1]{equation~(\protect\ref{eq:#1})}
\newcommand{\reqref}[1]{(\protect\ref{eq:#1})}
\newcommand{\rerem}[1]{Remark~\protect\ref{rem:#1}}

\newcommand{\sii}[2][1]{\ensuremath{\sigma_{#2}^{-{#1}}}}
\newcommand{\ssn}[2]{\ensuremath{\sigma_{#1}\cdots\sigma_{#2}}}
\newcommand{\aij}[2]{\ensuremath{B_{#1,#2}}}

\newcommand{\si}[2][{}]{\ensuremath{\sigma_{#2}^{#1}}}
\newcommand{\rh}[2][{}]{\ensuremath{\rho_{#2}^{#1}}}
\newcommand{\rhi}[2][1]{\ensuremath{\rho_{#2}^{-{#1}}}}
\newcommand{\rhs}[2]{\ensuremath{\rho_{#1}\cdots\rho_{#2}}}

\newcommand{\ssni}[2]{\ensuremath{\sigma_{#1}^{-1}\cdots\sigma_{#2}^{-1}}}

\begin{document}

\title{The braid groups of the projective plane}                    
\authors{Daciberg~Lima~Gon\c{c}alves\\John~Guaschi}                  
\coverauthors{Daciberg Lima Gon\noexpand\c{c}alves\\John Guaschi}
\asciiauthors{Daciberg Lima Goncalves and John Guaschi}

\address{Departamento de Matem\'atica - IME-USP\\Caixa 
Postal 66281 - Ag.~Cidade de S\~ao Paulo\\
CEP: 05311-970 - S\~ao Paulo - SP - Brasil\\{\rm 
and}\\Laboratoire de Math\'ematiques Emile Picard, UMR CNRS 5580 
UFR-MIG\\Universit\'e Toulouse III, 118, route de Narbonne\\
31062 Toulouse Cedex 4, France}

\asciiaddress{Departamento de Matematica - IME-USP\\Caixa 
Postal 66281 - Ag. Cidade de Sao Paulo\\
CEP: 05311-970 - Sao Paulo - SP - Brasil\\and\\Laboratoire de
Mathematiques Emile Picard, UMR CNRS 5580 
UFR-MIG\\Universite Toulouse III, 118, route de Narbonne\\
31062 Toulouse Cedex 4, France}

\gtemail{\mailto{dlgoncal@ime.usp.br}\qua{\rm 
and}\qua\mailto{guaschi@picard.ups-tlse.fr}}                     
\asciiemail{dlgoncal@ime.usp.br, guaschi@picard.ups-tlse.fr}

\begin{abstract}   % type your abstract below
Let $B_n(\rp)$ (respectively $P_n(\rp)$) denote the braid group
(respectively pure braid group) on $n$ strings of the real projective
plane $\rp$. In this paper we study these braid groups, in particular
the associated pure braid group short exact sequence of Fadell and
Neuwirth, their torsion elements and the roots of the `full twist'
braid. Our main results may be summarised as follows: first, the pure
braid group short exact sequence $1 \to
P_{m-n}(\rp\setminus\{x_1,\ldots, x_n\}) \to P_m(\rp) \to P_n(\rp) \to 1$ does
not split if $m\geq 4$ and $n=2,3$. Now let $n\geq 2$. Then in
$B_n(\rp)$, there is a $k$-torsion element if and only if $k$ divides
either $4n$ or $4(n-1)$. Finally, the full twist braid has a $k\th$
root if and only if $k$ divides either $2n$ or $2(n-1)$.
\end{abstract}

\asciiabstract{%
Let B_n(RP^2)$ (respectively P_n(RP^2)) denote the braid group
(respectively pure braid group) on n strings of the real projective
plane RP^2.  In this paper we study these braid groups, in particular
the associated pure braid group short exact sequence of Fadell and
Neuwirth, their torsion elements and the roots of the `full twist'
braid. Our main results may be summarised as follows: first, the pure
braid group short exact sequence 

1 --> P_{m-n}(RP^2 - {x_1,...,x_n}) --> P_m(RP^2) --> P_n(RP^2) --> 1 

does not split if m > 3 and n=2,3.  Now let n > 1.  Then in
B_n(RP^2), there is a k-torsion element if and only if k divides
either 4n or 4(n-1).  Finally, the full twist braid has a k-th
root if and only if k divides either 2n or 2(n-1).}

\primaryclass{Primary: 20F36, 55R80}                
\secondaryclass{Secondary: 55Q52, 20F05}              
\keywords{Braid group, configuration space, torsion, Fadell-Neuwirth short exact
sequence}

\maketitle 

%%%%%%%%%%%%%%%%%%%%   Start of main body of article

\section{Introduction}\label{sec:intro}

Braid groups of the plane $\Et$ were defined by Artin
in~1925~\cite{A1}, and further studied in~\cite{A2,A3}. They were
later generalised using the following definition due to
Fox~\cite{FoN}. Let $M$ be a compact, connected surface, and let
$n\in\N$. We denote the set of all ordered $n$-tuples of distinct
points of $M$, known as the \emph{$n\th$ configuration space of $M$},
by:
\begin{equation*}
F_n(M)=\{ (p_1,\ldots,p_n)\;\vert\; \text{$p_i\in M$ and $p_i\neq p_j$ if $i\neq j$} \}.
\end{equation*}
The configuration spaces $F_n(M)$ play an important r\^ole in several branches of mathematics and have been extensively studied, see~\cite{CG} for example. 

The symmetric group $\sn$ on $n$ letters acts freely on $F_n(M)$ by
permuting coordinates. The corresponding quotient will be denoted by
$D_n(M)$. Notice that $F_n(M)$ is a regular covering of $D_n(M)$. The
\emph{$n\th$ pure braid group $P_n(M)$} (respectively the
\emph{$n\th$ braid group $B_n(M)$}) is defined to be the fundamental
group of $F_n(M)$ (respectively of $D_n(M)$). If $m>n$ are positive
integers, then we may define a homomorphism \(
\map{p_{\ast}}{\pnm[m]}{\pnm} \) induced by the projection \(
\map{p}{F_m(M)}{F_n(M)} \) defined by \( p((x_1,\ldots, x_m))=
(x_1,\ldots,x_n) \). Representing $\pnm[m]$ geometrically
as a collection of $m$ strings, $p_{\ast}$ corresponds to forgetting
the last $(m-n)$ strings. 

\emph{We adopt the convention that,
unless explicitly stated, all homomorphisms\break \( \pnm[m] \to \pnm \)
in the text will be this one}.

If $M$ is without boundary, Fadell and Neuwirth study the map $p$, and show
\cite[Theorem~3]{FaN} that it is a locally trivial fibration. The
fibre over a point $(x_1,\ldots,x_n)$ of the base space is
$F_{m-n}(M\setminus\brak{x_1,\ldots,x_n})$ which we consider to
be a subspace of the total space via the map $\map
i{F_{m-n}(M\setminus\brak{x_1,\ldots,x_n})}{F_m(M)}$ defined by
$i((y_1,\ldots,y_{m-n}))=(x_1,\ldots,x_n,y_1,\ldots,y_{m-n})$.
Applying the associated long exact homotopy sequence, we obtain the
\emph{generalised pure braid group short exact sequence of Fadell and
Neuwirth}:
\begin{equation}\label{eq:split}
1 \to P_{m-n}(M\setminus\brak{x_1,\ldots,x_n}) \stackrel{i_{\ast}}{\to} \pnm[m] \stackrel{p_{\ast}}{\to}
\pnm \to 1, \tag{\normalfont\textbf{PBS}}
\end{equation}
where $n\geq 3$ if $M=\St$~\cite{Fa,FvB}, $n\geq 2$ if
$M=\rp$~\cite{vB}, and $n\geq 1$ otherwise~\cite{FaN}, and where
$i_{\ast}$ and $p_{\ast}$ are the homomorphisms induced by the maps
$i$ and $p$ respectively. The sequence also exists in the case  where $M$ is the $2$-disc $\dt$.

This short exact sequence has been widely studied, notably in recent
work in relation to the structure of the mapping class
groups~\cite{PR}. An interesting question is that of whether the
sequence~\reqref{split} splits. If $M$ is a $K(\pi,1)$ then the
existence of a section for $p_{\ast}$ is equivalent to that of a
cross-section for the fibration. In~\cite{A2}, Artin showed that if
$M=\dt$ then the sequence~\reqref{split} splits for all $n\in\N$. He
used this result to solve the word problem in $B_n(\dt)$ by showing
that $P_n(\dt)$ may be expressed as a repeated semi-direct product of
free groups.  Since then, the splitting problem has been studied for
other surfaces besides the plane. In the case $M=\St$, it was known
that there exists a section on the geometric level for all $n\geq
3$~\cite{Fa}. If $M$ is the $2$-torus then for all $n\geq 1$,\break Birman
exhibits an explicit algebraic section for the
sequence~\reqref{split}~\cite{Bi1}. However, for compact orientable
surfaces without boundary of genus $g\geq 2$, she poses the question
of whether the short exact sequence~\reqref{split} splits. The
authors of this paper have recently resolved this question, the
answer being positive if and only if $n=1$~\cite{GoG}. We will treat the case of non-orientable surfaces of genus at least two in a forthcoming paper~\cite{GG}. 

In this paper we shall study the braid groups of $\rp$, in particular
the splitting of the sequence~\reqref{split}, and the existence of a
section for the fibration $p$. These groups were first studied by Van
Buskirk~\cite{vB}, and more recently by Wang~\cite{Wa1}. Van Buskirk
showed that $P_1(\rp)=B_1(\rp)\cong \Z_2$, $P_2(\rp)$ is isomorphic
to the quaternion group $Q_8$, $B_2(\rp)$ is a dicyclic group of
order~$16$, and for $n>2$, $P_n(\rp)$ and $B_n(\rp)$ are infinite. He
also proved that these groups have elements of finite order
(including one of order~$2n$), such elements later being characterised
by Murasugi~\cite{M}. With respect to the splitting problem, Van
Buskirk showed that for $n\geq 2$, neither the fibration 
$\map{p}{F_n(\rp)}{F_1(\rp)}$ nor the homomorphism
$\map{p_{\ast}}{P_n(\rp)}{P_1(\rp)}$  admit a cross-section, but that the
fibration $\map{p}{F_3(\rp)}{F_2(\rp)}$ admits a cross-section (and
hence so does the corresponding homomorphism). In order to understand
better the configuration space $F_n(\rp)$, $n\in\N$, in
\repr{homotopy} we determine the homotopy type of its universal
covering space. From this, we may deduce its $i\th$ homotopy groups,
$i\geq 2$ (see \resec{rp2}):

\begin{prop}\label{prop:homorp2}
Let $n\in\N$, and let $i\geq 2$.
\begin{enumerate}
\item If $n=1$ then $\pi_i(F_1(\rp))\cong \pi_i(\St)$.
\item If $n\geq 2$ then $\pi_i(F_n(\rp))\cong \pi_i(\St[3])$.
\end{enumerate}
\end{prop}
We then show that  for $m\geq 4$ and $n=2,3$, neither the fibration
nor the short exact sequence~\reqref{split} admit a section. More
precisely:

\begin{thm}\label{th:nofmf2}
Let $m\geq 4$. Then:
\begin{enumerate}
\item\label{it:fmf2} the fibration $\map {p}{F_m(\rp)}{F_2(\rp)}$ does not admit a cross-section;
\item\label{it:fmf3} the fibration $\map {p}{F_m(\rp)}{F_3(\rp)}$ does not admit a cross-section. 
\end{enumerate}  
\end{thm}

\begin{thm}\label{th:nopnp2}
For $n=2,3$ and for all $m\geq 4$, the Fadell-Neuwirth pure braid group short exact sequence :
\begin{equation*}
1 \to P_{m-n}(\rp\setminus\brak{x_1,\ldots,x_n}) \stackrel{i_{\ast}}{\to} P_m(\rp) \stackrel{p_{\ast}}{\to} P_n(\rp) \to 1
\end{equation*}
does not split.
\end{thm}

Also in \resec{rp2}, we give an explicit algebraic section in the case
$m=3$ and $n=2$, which we use in \resec{rootpn} to study further the
structure of $P_3(\rp)$, notably its subgroups abstractly isomorphic
to $P_2(\rp)$.

Let $\dt\subseteq \rp$ be a topological disc. This inclusion induces
a (non-injective) homomorphism $\map {\iota}{B_n(\dt)}{B_n(\rp)}$.
Let $\sigma_1,\ldots,\sigma_{n-1}$ be a standard set of generators of
$B_n(\dt)$, and let:
\begin{equation}\label{eq:fulltwist}
\ft=(\sigma_1\cdots\sigma_{n-1})^n,
\end{equation}
be the \emph{full twist} braid of $B_n(\dt)$. By convention, we set $\ft[1]=1$. The braid $\iota(\ft)$
which for $n\geq 2$ is of order~$2$ and generates the centre $Z(B_n(\rp))$ of $B_n(\rp)$~\cite{vB},
shall also be denoted by $\ft$, and will be called \emph{the full
twist of $B_n(\rp)$}. If $k\in\N$, then we shall say that
$\alpha\in B_n(M)$ is a \emph{$k\th$ root of $\ft$} if $\alpha^k=\ft$
and  $\alpha^j\neq\ft$ for all $1\leq j<k$. As we previously
mentioned, the torsion elements of $B_n(\rp)$, as well as the roots
of the full twist were characterised by Murasugi. However, in both
cases, the orders of the elements are not completely clear, even for
$P_n(\rp)$. In \resec{rootpn}, we study the torsion as well as the
roots of $\ft$ in $P_n(\rp)$ using methods different to those of
Murasugi, and prove that for all $n\geq 2$ the torsion of $P_n(\rp)$
is precisely~$2$ and~$4$ (see \reco{torspn}). In \resec{rootbn}, we
do the same for $B_n(\rp)$:

\begin{thm}\label{th:torsbn}
Let $n\geq 2$. Then $B_n(\rp)$ has an element of order $\ell$ if and 
only if $\ell$ divides either $4n$ or $4(n-1)$.
\end{thm}

From \repr{unibn} (see also \resec{rootbn}), $\ft$ is the unique
element of $B_n(\rp)$ of order~$2$. It thus follows that:

\begin{thm}\label{th:rootbn}
Let $n\geq 2$. Then the full twist $\ft$ has a $k\th$ root in
$B_n(\rp)$ if and only if $k$ divides either $2n$ or $2(n-1)$.
\end{thm}

\begin{ack}
This work began  during the visit of the second author to the
Departmento de Matem\'atica do  IME-Universidade de S\~ao Paulo
during the period 7\up{th}~July~--~3\up{rd}~August 2002.   The visit
was supported by the international Cooperation Capes/Cofecub project
number 364/01.
\end{ack}

\section{The splitting problem for $\rp$}\label{sec:rp2}

In this section, we recall, reformulate and generalise some results
obtained in~\cite{Fa,vB} concerning the configuration spaces
$F_n(\rp)$ of the projective plane. If $n\geq 2$, we show that the
universal covering of the configuration spaces $F_n(\rp)$ has the
homotopy type of the $3$-sphere $\St[3]$. We then prove that if
$n\geq 2$, the problem of finding a cross-section for the fibrations
$F_m(\rp)\stackrel{p}{\to} F_n(\rp)$ introduced in \resec{intro} is
equivalent to that of the splitting of the corresponding pure braid
group short exact sequence. Finally, we prove
Theorems~\ref{th:nofmf2} and~\ref{th:nopnp2}, that is, if $n=2,3$ and
$m\geq 4$ then the fibration $F_m(\rp)\stackrel{p}{\to} F_n(\rp)$
does not admit a cross-section, with a similar conclusion for the
splitting of the short exact sequence~\reqref{split}.

\subsection{The configuration spaces $F_n(\rp)$ and a description of their homotopy groups}

We first determine the homotopy type of the universal covering of the configuration spaces 
$F_n(\rp)$ and discuss their fundamental groups. Partial results in this direction were obtained 
in~\cite{Fa,vB}.

\begin{prop}\mbox{}\label{prop:homotopy}
\begin{enumerate}
\item \label{it:proiso} Let $m>n\geq 2$. Then for each $i>1$, the
projection $\map p{F_m(\rp)}{F_n(\rp)}$ onto the first $n$
coordinates induces an isomorphism of the corresponding $i\th$
homotopy groups.
\item \label{it:f1rp} For $n=1$, the universal covering of $F_n(\rp)$ is the sphere $\St$. The fundamental group of $F_1(\rp)$ is $\Z_2$.
\item \label{it:f2rp} For $n=2$, the fundamental group of $F_n(\rp)$
is isomorphic to the quaternion group $Q_8$ of order~$8$.
\item \label{it:f3rp} For $n\geq 3$, the fundamental group of $F_n(\rp)$ is infinite.
\item \label{it:fnrp} For $n\geq 2$, the homotopy type of the
universal covering of $F_n(\rp)$ is that of the $3$-sphere $\St[3]$. 
\item \label{it:orbit} The orbit space arising from the standard
action of $Q_8$ on $\St[3]$ has the same homotopy type as $F_2(\rp)$.
\end{enumerate}
\end{prop}

Notice that \repr{homotopy} implies \repr{homorp2}.

\begin{proof}[Proof of \repr{homotopy}]
To prove part~{(\ref{it:proiso})}, let $m>n\geq 2$. For $i>1$, consider the long exact 
sequence in homotopy associated to the (locally trivial) fibration $\map p{F_m(\rp)}{F_n(\rp)}$ 
given by $p((x_1,\ldots,x_m))=(x_1,\ldots,x_n)$:
\begin{gather*}
\cdots \to \pi_i(F)\to \pi_i(F_m(\rp)) \to \pi_i(F_n(\rp)) \to \pi_{i-1}(F) \to\cdots\\
\cdots \to \pi_2(F)\to \pi_2(F_m(\rp)) \to \pi_2(F_n(\rp)) \to \pi_1(F)\to\cdots
\end{gather*}
where $F=F_{m-n}(\rp\setminus\brak{x_1,\dots,x_n})$ is the fibre over
$(x_1,\dots,x_n)\in F_n(\rp)$. Since $F$ is an aspherical
space~\cite{FaN}, and $\pi_2(F_q(\rp))=0$ for all $q\geq
2$~\cite{vB}, it follows that $p$ induces an isomorphism of
$\pi_i(F_m(\rp))$ onto $\pi_i(F_n(\rp))$.

part~{(\ref{it:f1rp})} is clear, since $F_1(\rp)=\rp$, and parts~{(\ref{it:f2rp})} and~{(\ref{it:f3rp})} follow from~\cite{vB}. 

Before proving part~{(\ref{it:fnrp})}, let us recall the
following general facts. For $m\geq 1$, $F_m(\rp)$ is a
$2m$-dimensional smooth manifold~\cite{CJ}, and so there exists a
CW-complex $K_m$ (of dimension less than or equal to $2m$) whose
homotopy type is that of $F_m(\rp)$~\cite{Hir}. Let $\widetilde{K_m}$
(respectively $\widetilde{F_m(\rp)}$) denote the universal covering
of $K_m$ (respectively $F_m(\rp)$). Then $\widetilde{K_m}$ is also a
CW-complex whose homotopy type is that of
$\widetilde{F_m(\rp)}$~\cite{Hir}. 

To prove part~{(\ref{it:fnrp})}, we first prove the statement
for $n=2$. Let $\left\{ \begin{alignedat}{2} & \St &\to &\rp\\ & v
&\longmapsto &[v] \end{alignedat}\right.$ denote the usual covering
map which identifies antipodal points. \vrule width 0pt height 12pt
Consider the map $\left\{
\begin{alignedat}{2} \varphi\colon & V_{3,2}&\to & F_2(\rp)\\ &
(e_1,e_2)& \longmapsto & ([e_1],[e_2]) \end{alignedat}\right.$ where
$V_{3,2}$ is the Stiefel manifold of orthogonal $2$-frames in
Euclidean $3$-space. Then $\map {\varphi}{V_{3,2}}{\varphi(V_{3,2})}$
is a covering map with discrete fibre consisting of $4$ points, and
so $V_{3,2}$ and $\varphi(V_{3,2})$ have the same homotopy type.
Further, $\varphi(V_{3,2})$ is a deformation retract of $F_2(\rp)$;
indeed, if $([v_1],[v_2])\in F_2(\rp)$ then one may obtain an element
of $\varphi(V_{3,2})$ by deforming the second vector $v_2$ within the
plane containing $v_1$ and $v_2$ to a vector $v_2'$ orthogonal to
$v_1$. Although $v_2'$ is not unique in $V_{3,2}$, its image in $\rp$
is, and this gives rise to a (well-defined) deformation retract of
$F_2(\rp)$ onto $\varphi(V_{3,2})$. It follows that $V_{3,2}$ and
$F_2(\rp)$ have the same homotopy type. So their universal covers do
too, and the result follows from the fact that $V_{3,2}$ is
homeomorphic to $\rp[3]$ whose universal covering is $\St[3]$.

Now suppose that $n\geq 3$. Consider the fibration
$\map{p}{F_n(\rp)}{F_2(\rp)}$. By part~{(\ref{it:proiso})}, $p$
induces an isomorphism of the $i\th$ homotopy groups of $F_n(\rp)$
and $F_2(\rp)$ for $i\geq 2$. Composing $p$ appropriately with
homotopy equivalences between $K_n$ and $F_n(\rp)$, and $K_2$ and
$F_2(\rp)$, it follows that there exists a map $K_n\to K_2$ which
induces an isomorphism of their $i\th$ homotopy groups for $i\geq 2$,
and thus a map $\widetilde{K_n}\to \widetilde{K_2}$ which induces an
isomorphism of their $i\th$ homotopy groups for $i\geq 2$. From
Whitehead's theorem, we see that $\widetilde{K_n}$ and
$\widetilde{K_2}$ are homotopy equivalent~\cite{Wh2}. Since
$\widetilde{K_n}$ and $\widetilde{F_n(\rp)}$ (respectively
$\widetilde{K_2}$ and $\widetilde{F_2(\rp)}$) have the same homotopy 
type, then so do $\widetilde{F_n(\rp)}$ and
$\widetilde{F_2(\rp)}$, and this completes the proof of
part~{(\ref{it:fnrp})}.

Finally, for part~{(\ref{it:orbit})}, consider the proof of
part~{(\ref{it:fnrp})} in the case $n=2$. The covering map $\map
{\varphi}{V_{3,2}}{\varphi(V_{3,2})}$ is induced by the natural free
action of $\Z_2\oplus \Z_2$ on $V_{3,2}$ given as follows: if
$\epsilon_i\in\brak{\pm 1}$ for $i=1,2$ then
$(\epsilon_1,\epsilon_2)\cdot
(e_1,e_2)=(\epsilon_1e_1,\epsilon_2e_2)$.  Recall that $\St[3]$ may
be considered as the topological group of quaternions of unit
modulus. One may consider the quotient of $\St[3]$ arising from the
standard action of $Q_8$ in two stages: first quotient $\St[3]$ by
the centre $\Z_2$ of $Q_8$ to obtain $\rp[3]$ which is homeomorphic
to $V_{3,2}$, then quotient $V_{3,2}$ by $Q_8/\Z_2\cong \Z_2\oplus
\Z_2$. We saw above that this quotient has the same homotopy type as
$F_2(\rp)$, and this completes the proof of
part~{(\ref{it:orbit})} and thus that of the proposition. 
\end{proof}   

\subsection{Generalities on the splitting problem}

We now consider the problem of determining a cross-section to the fibration $F_m(\rp) \stackrel{p}{\to} F_n(\rp)$ defined in \resec{intro}. If $n\geq 2$ then this problem is equivalent to showing that the corresponding pure braid group short exact sequence splits:

\begin{prop}\label{prop:xsplit}
Let $m>n\geq 2$. Then the fibration $\map p{F_m(\rp)}{F_n(\rp)}$
admits a cross-section if and only if the corresponding pure braid
group short exact sequence
\begin{equation*}
1 \to P_{m-n}(\rp\setminus\{x_1,\ldots,x_n\}) \to P_m(\rp) \stackrel{p_{\ast}}{\to} P_n(\rp) \to
1
\end{equation*}
splits.
\end{prop}

\begin{proof}
The `only if' part is clear. For the `if part', the existence of the
given short exact sequence and its splitting imply that we can find a
section over the $1$-skeleton, and then over the $2$-skeleton of
$F_n(\rp)$~\cite[Theorem~4.3.1]{Ba}. From~\cite{Wh1}, the
obstructions to extending the map over the $q$-skeleta, $q\geq 3$,
lie in cohomology classes whose coefficients belong to the $(q-1)\th$
homotopy group of the fibre $F_{m-n}(\rp\setminus\{x_1,\ldots,x_n\})$. But
this fibre is an Eilenberg-MacLane space.  Hence we can extend to a
section over the entire space, and the result follows.
\end{proof}

The cases not covered by \repr{xsplit} may be resolved as follows:
\begin{prop}\label{prop:nosplit}\mbox{}
\begin{enumerate}
\item \label{it:nosplit} Let $m\geq 2$. Then the surjective group homomorphism $\map{p_{\ast}}{P_m(\rp)}{P_1(\rp)}$ does not admit a section.
\item\label{it:noxm2} Let $m\geq 2$. Then the fibration $\map
p{F_m({\rp})}{F_1({\rp})}$ does not admit a cross-section.
\end{enumerate}
\end{prop}

\proof
\begin{enumerate}
\item We first note that by the commutativity of the following diagram,
\begin{equation*}
\begin{CD}
P_m(\rp) @>>> P_1(\rp)\\
@VVV @|  \\
P_2(\rp) @>>> P_1(\rp)
\end{CD}
\end{equation*}
the existence of a section for $P_m(\rp) \to P_1(\rp)$ implies that
there exists a section for $P_2(\rp) \to P_1(\rp)$
(cf~\cite[Proposition~3]{GoG}), and hence it suffices to
suppose that $m=2$. To see that $P_2(\rp) \to P_1(\rp)$ does not
admit a section, notice that $P_2(\rp)\cong Q_8$ has only one element
of order $2$, the full twist $\ft[2]$. But this element is sent to
the trivial element of $P_1(\rp)\cong\Z_2$ under the projection onto
$P_1(\rp)$,  and so there is no section. 
\item If there were to exist a cross-section for $F_m({\rp})\to
F_1({\rp})$ then there would exist a section for $P_m(\rp) \to
P_1(\rp)$, which contradicts part~{(\ref{it:nosplit})}.\qed
\end{enumerate}

\begin{rem}
The fact that the fibration given in part~{(\ref{it:noxm2})} of
\repr{nosplit} does not admit a cross-section also follows from the
fact that $\rp$ has the fixed point property (see also~\cite{vB}).
\end{rem}

\subsection{The splitting problem for $n\in\brak{2,3}$, and the structure of $P_3(\rp)$}\label{sec:n23}

This problem of the existence of a cross-section for the fibrations
$F_m(M)\to F_n(M)$ was considered by Fadell and Neuwirth during the
1960's~\cite{Fa,FaN}. The case $M=\rp$ was later studied by Van
Buskirk who constructed an explicit cross-section when $m=3$ and
$n=2$~\cite{vB}. He posed the question of the characterisation of the
values of $m$ and $n$ for which there exists a cross-section. As far as we
know, if $m\geq 4$ and $n\geq 2$, the question is still completely
open.

We first consider the case $m=3$ and $n=2$ of Van Buskirk, and
construct an explicit splitting of the corresponding pure braid group
short exact sequence. This will be useful later on when we study
$P_3(\rp)$ in more detail. We show that for $m>3$ and $n=2,3$, the
pure braid group short exact sequence does not split. We obtain this
result by showing that for $m>3$, the fibration $F_m(\rp) \to
F_2(\rp)$ does not admit a cross-section. 

The cross-section $\map{\sigma}{F_{2}(\rp)}{F_{3}(\rp)}$ defined by Van Buskirk to the
fibration $F_{3}(\rp) \to F_{2}(\rp)$ is
given  by $\sigma([v_1],[v_2])= ([v_1],[v_2],[v_1\bigwedge v_2])$, where
$\bigwedge$ is the vector product (cf the proof of
part~{(\ref{it:fnrp})} of \repr{homotopy}). Hence the following short
exact sequence splits:
\begin{equation}\label{eq:p3split}
1 \to \pi_1(\rp\setminus\{ x_1,x_2 \}) \to P_3(\rp) \stackrel{p_{\ast}}{\to} P_2(\rp) \to 1,
\end{equation}
where $p$ is induced by the projection $F_{3}(\rp) \to F_{2}(\rp)$
onto the last two coordinates. Since the kernel is isomorphic to
$\F$, the free group on two generators, it follows that $P_3(\rp)$ is
isomorphic to the semi-direct product of $\F$ with $Q_8$.

If $n\in\N$, a generating set for $P_n(\rp)=
\pi_1(F_n(\rp),(x_1,\ldots,x_n))$ may be obtained as follows. By
removing a disc from $\rp$, one obtains a disc with a twisted band
attached. For $1\leq i\leq n$, let $\rh{i}$ be a braid represented by
the loop based at the point $x_i$ which goes round the band, while
the other points $x_j$, $j\neq i$, stay fixed. Let $\dt\subseteq\rp$
be a topological disc which contains $\brak{x_1,\ldots,x_n}$.
This inclusion induces a (non-injective) homomorphism $B_n(\dt)\to
B_n(\rp)$. Let $\si{1},\ldots,\si{n-1}\in B_n(\rp)$ be the images of
the usual Artin generators of $B_n(\dt)$ under this homomorphism. For
$1\leq j<k\leq n$, set $\aij{j}{k} = \ssn{k-1}{j+1} \si[2]{j}
\ssni{j+1}{k-1}\in P_n(\rp)$. Then the union of the $\rh{i}$ and
$\aij{j}{k}$ generates $P_n(\rp)$ (cf~\cite{vB,S}).

From~\req{p3split}, we may obtain a presentation of $P_3(\rp)$, using
the fact that it is an extension~\cite{J}. It has generators 
$\brak{\rh{1},\rh{2},\rh{3},B_{2,3}}$ where $\brak{\rh{1},\rh{2}}$
are (coset representatives of) generators of the quotient $P_2(\rp)$,
and $\{\rh{3},B_{2,3}\}$ generates the kernel. A set of relations is
given by:
\begin{equation}\label{eq:relp3}
\left.
\begin{aligned}
\rh{1}B_{2,3}\rhi{1}&=B_{2,3} \quad & \rh{2}B_{2,3}\rhi{2}&=\rhi{3}B_{2,3}^{-1}\rh{3}\\
\rhi{1}B_{2,3}\rh{1}&=B_{2,3} \quad & \rhi{2}B_{2,3}\rh{2}&=
B_{2,3}^{-1}\rh{3} B_{2,3}\rhi{3}B_{2,3}\\
\rh{1} \rh{3}\rhi{1}&=\rhi{3}B_{2,3}  \quad & 
\rh{2} \rh{3}\rhi{2}&=\rhi{3}B_{2,3}^{-1}\rh[2]{3}\\
\rhi{1} \rh{3}\rh{1}&=B_{2,3}\rhi{3}  \quad & 
\rhi{2} \rh{3}\rh{2}&=B_{2,3}^{-1}\rh{3}\\
\rh{1}\rh{2}\rh{1}\rhi{2} &=B_{2,3}^{-1}\rh[2]{3} \quad &   
\rhi[2]{1}\rh[2]{2} &= B_{2,3}\rhi[2]{3} B_{2,3}.
\end{aligned}
\right\}
\end{equation}
This presentation was obtained using the following useful relations: 
\begin{align*}
\rh[2]{1} &=B_{1,2}B_{1,3}, &
\rh[2]{2}&=B_{1,2}B_{2,3}, & \rh[2]{3}&=B_{1,3}B_{2,3},
\end{align*}
and the presentation $\ang{\rh{1},\rh{2}\; \left\vert \; \rh[2]{1}= \rh[2]{2}, \, \rh{2}\rh{1}\rhi{2}= \rhi{1}\right.}$ of $B_2(\rp)$.

Now we will define an algebraic section. It will be also be used in
the following section to study properties of the torsion elements
of $P_3(\rp)$. 

\begin{prop}\label{prop:p2p3}
The function $\map s{\{\rh{1},\rh{2}\}}{P_3(\rp)}$ defined by
$s(\rh{1})=\rh{1}\rh{3}$ and $s(\rh{2})= \rh{3}\rh{2}$ extends to a
homomorphism $\map s{P_2(\rp)}{P_3(\rp)}$, and gives a splitting of
the short exact sequence~\reqref{p3split}.
\end{prop}

\begin{proof} 
The set $\brak{\rh{1},\rh{2}}$ generates the group $P_2(\rp)$. One
may check that $s(\rh[2]{1})=s(\rh[2]{2})$ and
$s(\rh{2}\rh{1}\rhi{2})= s(\rhi{1})$. Since $P_2(\rp)\cong Q_8$, this
implies that $s$ is a homomorphism which composed with the
homomorphism $p_{\ast}$ of the short exact sequence~\reqref{p3split}
gives the identity. 
\end{proof}

\begin{cor}\label{cor:repres}
The representation $\map{\phi}{Q_8}{\aut{\F}}$ induced by the
section $s$ of \repr{p2p3} is given by $\phi(\rh{1})=\chi_1$,
$\phi(\rh{2})=\chi_2$, and $\phi(\rh{1}\rh{2})=\chi_3$, where:
\begin{align*}
\chi_1(x)&=y & \chi_1(y)&=x,\\
\chi_2(x)&=y^{-1} & \chi_2(y)&=x^{-1},\\
 \chi_3(x)&=x^{-1} & \chi_3(y)&=y^{-1},
\end{align*}
and where $x=\rh{3}$ and $y=\rhi{3}B_{2,3}$ form a basis of $\F$. 
\end{cor}

\begin{proof}
From above, we know that $\brak{\rhi{3},B_{2,3}}$ is a basis of
$\F$.  The homomorphism $s$ induces a representation $\phi$ of
$Q_8$ in the group of automorphisms of $\F$. For $i=1,2$, let
$\chi_i$ denote the automorphism $\phi(\rh{i})$. Then we obtain:
\begin{align*}
\chi_1(\rh{3})&=\rhi{3}B_{2,3}, & \chi_1(B_{2,3})&=\rhi{3}B_{2,3}\rh{3}, &   \chi_2(\rh{3})&=B_{2,3}^{-1}\rh{3}, & \chi_2(B_{2,3})&=B_{2,3}^{-1}.
\end{align*}    
Consider the basis $x=\rh{3}$, $y=\rhi{3}B_{2,3}$ of $\F$.
With respect to $\brak{x,y}$, the automorphisms $\chi_1$,
$\chi_2$ and $\chi_3=\chi_2\circ \chi_1$ are as given in the
statement of the corollary. 
\end{proof}

As a preliminary step to proving \reth{nofmf2}, let us show that
$P_4(\rp) \stackrel{p_{\ast}}{\to} P_2(\rp)$ does not split.

\begin{prop}\label{prop:coinfree}
Let $\map{f_1,f_2}{F_{2}(\rp)}{\rp}$ be a pair of coincidence-free
maps such that the induced homomorphism
$\map{f_{1\ast}}{P_2(\rp)}{P_1(\rp)}$ is surjective. Then
$\map{f_{2\ast}}{P_2(\rp)}{P_1(\rp)}$ is also surjective, and
$f_{1\ast}\neq f_{2\ast}$ ie, the induced homomorphisms on the
fundamental group are different and both non trivial.
\end{prop}

\begin{proof}
Consider the map $\left\{ \begin{alignedat}{3} (f_1,f_2)\colon
F_{2}(\rp) &\to \rp \times \rp\\ w & \longmapsto  (f_1(w),f_2(w)).
\end{alignedat}\right.$ If $(f_1,f_2)$ is coincidence free (or can be
deformed to a coincidence-free pair) then there is a  factorisation
via $F_2(\rp)=\rp\times\rp\setminus\Delta$, where $\Delta$ is the
diagonal, and so there exists a homomorphism
$\map{\phi}{P_2(\rp)}{P_2(\rp)}$ which makes the following diagram
commute: 
\begin{equation*}
\xymatrix{%
& & \ar[d]^{g_{\ast}} P_{2}(\rp)\\
P_{2}(\rp) \ar@{-->}[urr]+DL^(.6){\phi}  \ar[rr]_(.38){(f_1,f_2)_{\ast}} & & \pi_1(\rp) \times \pi_1(\rp),
}
\end{equation*}
where $g_{\ast}$ is the homomorphism induced by the inclusion
$\map[\lhra] g{F_2(\rp)}{\rp\times\rp}$. Identifying $\pi_1(\rp
\times \rp)$ with $\Z_2+\Z_2$ in the obvious way, let $\rh{1},\rh{2}$
be generators of $P_2(\rp)$ such that $g_{\ast}(\rh{1})=(1,0)$ and 
$g_{\ast}(\rh{2})=(0,1)$. Let us first show that $f_{2\ast}$ cannot
be the trivial homomorphism. 
Suppose on the contrary that $f_{2\ast}=0$. If $a\in P_2(\rp)$ then
there exist $i\in\brak{0,1,2,3}$ and $j\in\brak{0,1}$ such that
$\phi(a)=\rh[i]{1}\rh[j]{2}$. Since $g_{\ast}\circ
\phi=(f_1,f_2)_{\ast}$, it follows that $j=0$. Now $f_{1\ast}$ is
surjective, so there exist $b\in\brak{\rh{1},\rh{2}}$ and
$\epsilon\in\brak{\pm 1}$ such that $\phi(b)=\rh[\epsilon]{1}$. Hence
$P_2(\rp)/\ker{\phi}\cong\Z_4$. But $P_2(\rp)\cong Q_8$ of which the
only quotient of order~$4$ is $\Z_2+\Z_2$~--~a contradiction. Since
$P_1(\rp)\cong\Z_2$, we conclude that $f_{2\ast}$ is surjective.

Suppose now that $f_{1\ast}= f_{2\ast}$. Let $a\in P_2(\rp)$. Then
there exist $i\in\brak{0,1,2,3}$ and $j\in\brak{0,1}$ such that
$\phi(a)=\rh[i]{1} \rh[j]{2}$. Since $g_{\ast}\circ
\phi=(f_1,f_2)_{\ast}$, it follows that $i \equiv j\bmod 2$. Now
$f_{1\ast}$ is surjective, so there exist $b\in\brak{\rh{1},\rh{2}}$
and $\epsilon\in\brak{\pm 1}$ such that
$\phi(b)=\rh[\epsilon]{1}\rh{2}$. By an argument similar to that of
the previous paragraph, we obtain a contradiction, whence
$f_{1\ast}\neq f_{2\ast}$.
\end{proof}

We are now able to prove the main results of this section.

\begin{proof}[Proof of \reth{nofmf2}(\ref{it:fmf2})]
As in the proof of part~{(\ref{it:noxm2})} of \repr{nosplit}, it
suffices to show that there is no section for $m=4$. Suppose on the
contrary that such a section exists. Then there exists a map
$\map{\xi}{F_2(\rp)}{F_4(\rp)}$ of the form 
$\xi(x_1,x_2)=(x_1,x_2,f_1(x_1,x_2),f_2(x_1,x_2))$.  For $i=1,2$, let
$\map{p_i}{F_2(\rp)}{\rp}$ denote the map $p_i(x_1,x_2)=x_i$.  Then
the four maps $\map{f_1,f_2,p_1,p_2}{F_2(\rp)}{\rp}$ are pairwise
coincidence free, and so by \repr{coinfree}, the induced
homomorphisms are all surjective and pairwise distinct. But, since
$P_2(\rp)\cong Q_8$ and $P_1(\rp)\cong \Z_2$, it follows that there
are only three surjective homomorphisms $P_2(\rp)\to P_1(\rp)$. We
thus obtain a contradiction, and the result follows.
\end{proof}

\begin{proof}[Proof of \reth{nopnp2}]
First let $n=2$. By \repr{xsplit}, if this sequence were to split then there would exist a cross-section to the given fibration. But this would contradict \reth{nofmf2}{(\ref{it:fmf2})}.

Now let $n=3$. Consider the following commutative diagram:
\begin{equation*}
\begin{CD}
P_m(\rp) @>>> P_3(\rp)\\
@| @VVV \\
P_m(\rp) @>>> P_2(\rp). 
\end{CD}
\end{equation*}
Since the right-hand vertical arrow admits a section, if there were to exist a section for $P_m(\rp)\to P_3(\rp)$ then there would also exist a section for $P_m(\rp)\to P_2(\rp)$ which contradicts the previous paragraph.
\end{proof}  

\begin{proof}[Proof of \reth{nofmf2}(\ref{it:fmf3})]
The result follows from \repr{xsplit} and \reth{nopnp2}. This also completes the proof of \reth{nofmf2}.
\end{proof}

\section{Torsion elements and roots of the full twist in $P_n(\rp)$}\label{sec:rootpn}

Let $\ft$ denote the full twist braid of $B_n(\rp)$ as defined in
\req{fulltwist}. In~\cite{M},  Murasugi gave a characterisation of
the torsion elements of $B_n(\rp)$ and the quotient of $B_n(\rp)$ by its centre
$\ang{\ft}$. However, the orders of the elements are
not clear, even in the case of $P_n(\rp)$. In this section, we study
the torsion and the roots of $\ft$ in $P_n(\rp)$ with methods
different to those of Murasugi. Using the pure braid group short
exact sequence, we first show that the possible torsion of $P_n(\rp)$
is~$2$ or~$4$. We then go on to prove that $\ft$ is the unique
element of $P_n(\rp)$ of order~$2$, and that for all $n\geq 2$,
$P_n(\rp)$ possesses elements of order~$4$. It will follow that the
full twist possesses $k\th$ roots in $P_n(\rp)$, $k\geq 2$, if and
only if $k=2$. Finally, using the semi-direct product decomposition
given in \repr{p2p3}, we study in more detail the torsion elements of
$P_3(\rp)$.

We first observe that the possible torsion of $P_n(\rp)$ is either~$2$ or~$4$. 

\begin{lem}\label{lem:notors}
Let  $1 \to K \to G \stackrel{f}{\to} H \to 1$ be a short exact
sequence of groups such that $K$ is torsion free. If $x\in G$ is a
torsion element of order $l$, then the element $f(x)$ of $H$ has also
order $l$.
\end{lem}

\begin{proof}
Let $x\in G$ be a torsion element of order $l$. Therefore $f(x)^l=1$,
and the order $k$ of $f(x)$ divides $l$. Consider $x^k$. It belongs
to $\ker{f}=K$, and $(x^k)^{\frac{l}{k}}=1$. Since $K$ is torsion
free, it follows that $x^k=1$, and we conclude that $k=l$.
\end{proof}

\begin{cor}\label{cor:tor24}
Let $n\in\N$. Then the possible (non-trivial) torsion of $P_n(\rp)$ is~$2$ or~$4$.
\end{cor}

\begin{proof}
Since $P_1(\rp)\cong\Z_2$ and $P_2(\rp)\cong Q_8$, the result follows
if $n\in\brak{1,2}$. If $n\geq 3$,  consider the following pure braid
group short exact sequence:
\begin{equation*}
1\to P_{n-2}(\rp\setminus\brak{x_1,x_2}) \to P_n(\rp) \stackrel{p_{\ast}}{\to} P_2(\rp)\to 1. 
\end{equation*}
Now $F_{n-2}(\rp\setminus\brak{x_1,x_2})$ is a finite-dimensional
Eilenberg-MacLane space, and so by a theorem of P.~Smith,
$P_{n-2}(\rp\setminus\brak{x_1,x_2})$ is torsion free~\cite{FaN}. The
result follows from \relem{notors}.
\end{proof}

\begin{prop}\label{prop:fullpn}
Let $n\geq 2$. Then the full twist $\ft$ is the unique element of
$P_n(\rp)$ of order~$2$ in $P_n(\rp)$. In particular, the square of any
element of $P_n(\rp)$ of order~$4$ is the full twist.
\end{prop}

\begin{proof}
The proof of the first part is by induction. The result is true for
$n=2$, since $P_2(\rp)\cong Q_8$, and the unique element of order~$2$ of this group is the full twist. Suppose that the result is true for some $n\geq 2$. Consider the
following short exact sequence:
\begin{equation*}
1 \to F \to P_{n+1}(\rp) \stackrel{p_{\ast}}{\to} P_n(\rp) \to 1,
\end{equation*}
where $F=\pi_1(\rp\setminus\brak{x_1,\ldots,x_n})$ is torsion free.
If $\alpha\in P_{n+1}(\rp)$ is of order~$2$ then $p_{\ast}(\alpha)\in
P_n(\rp)$ is also of order~$2$ by \relem{notors}, and so
$p_{\ast}(\alpha)=\ft$ by the induction hypothesis. Since
$p_{\ast}(\ft[n+1])=\ft$, it follows that $\alpha=\ft[n+1]\cdot w$,
where $w\in F$. Then $1=\alpha^2=\ft[n+1]^2\cdot w^2=w^2$ since
$\ft[n+1]$ generates $Z(P_{n+1}(\rp))$ and is of
order~$2$~\cite{vB,M}. But $F$ is torsion free, thus $w=1$, and so
$\ft[n+1]$ is the unique element of $P_{n+1}(\rp)$ of order~$2$. The
fact that the square of any element of $P_n(\rp)$ of order~$4$ is the full twist
then follows easily.
\end{proof}  

Since $P_2(\rp)$ has an element of order~$4$, it follows from
\repr{p2p3} that $P_3(\rp)$ does too. In spite of the fact that 
$P_n(\rp) \stackrel{p_{\ast}}{\to} P_2(\rp)$ does not admit a
section, we will show that  $P_n(\rp)$ contains elements of order~$4$
for all $n\geq 2$. We first recall a presentation of $B_n(\rp)$.

\begin{prop}[Van Buskirk~\cite{vB}]\label{prop:present}
The following constitutes a presentation of the group $B_n(\rp)$:

{\bf Generators}\vspace{-10pt}
$$\si{1},\ldots,\si{n-1},\rh{1},\ldots,\rh{n}.$$\eject

{\bf Relations}\vspace{-10pt}
\begin{align*}
\si{i}\si{j} &=\si{j}\si{i}\quad\text{if $\lvert i-j\rvert\geq 2$,}\\
\si{i}\si{i+1}\si{i}&=\si{i+1}\si{i}\si{i+1} \quad\text{for $1\leq i\leq n-2$,}\\
\si{i}\rh{j}&=\rh{j}\si{i}\quad\text{for $j\neq i,i+1$,}\\
\rh{i+1}&=\sii{i}\rh{i}\sii{i} \quad\text{for $1\leq i\leq n-1$,}\\
\rh[-1]{i+1}\rh[-1]{i}\rh{i+1}\rh{i}&= \si[2]{i} \quad\text{for $1\leq i\leq n-1$,}\\
\rh[2]{1}&=\si{1}\si{2}\cdots\si{n-2}\si[2]{n-1} \si{n-2}\cdots\si{2}\si{1}.
\end{align*}
\end{prop}

We warn the reader that the $\rh{i}$ of \repr{present} are slightly different from the $\rh{i}$ defined in \resec{n23}. From this presentation, if $2\leq j\leq n$, we obtain the
relation 
\begin{equation}\label{eq:rhj1}
\rh{j}=\ssni{j-1}{1}\rh{1}\ssni{1}{j-1}.
\end{equation}
Notice also that $\aij{i}{j+1}= \si{j}\aij{i}{j} \sii{j}$ for all $1\leq i<j\leq n-1$.

\begin{lem}\label{lem:commute}
Let $1\leq i<j\leq n$. Then $\rh{j}\rh{i}=\rh{i}\rh{j}\aij{i}{j}$.
\end{lem}

\begin{proof}
Let $1\leq i\leq n-1$. The proof is by induction on $j-i$. Suppose
first that $j-i=1$. Then $\rh{i+1}\rh{i}=\rh{i}\rh{i+1}\si[2]{i}
=\rh{i}\rh{i+1}\aij{i}{i+1}$ by \repr{present}. Suppose now that the
given equation holds for $i$ and $j$ satisfying $1\leq i<j\leq n-1$
and $1\leq j-i\leq n-2$. Applying \repr{present}, we have:
\begin{align*}
\rh{j+1}\rh{i} &=  \sii{j}\rh{j}\sii{j}\rh{i}=\sii{j}\rh{j}\rh{i}\sii{j} =\sii{j}\rh{i}\rh{j}\aij ij\sii{j}\\
&=
\rh{i}\sii{j}\rh{j}\sii{j}\si{j}\aij ij\sii{j} =  
 \rh{i}\rh{j+1}\aij{i}{j+1},
\end{align*}
which proves the lemma.
\end{proof}

\begin{prop}\label{prop:rhosq}
Let $n\in\N$. In $B_n(\rp)$, $(\rhs {n-1}1)^2=\ft$.
\end{prop}

\begin{proof}
If $n=1$ then both sides are trivial, and the result is true. So let $n\geq 2$.
Let us first show that $(\rhs {n-1}1)^2$ is central in $B_n(\rp)$. From \repr{present}, it suffices to prove that it commutes with each member of the generating set $\brak{\si{1},\ldots, \si{n-1},\rh{n}}$. By \relem{commute} we have that $\si{i}\rh{i+1}\rh{i}\sii{i} = \si{i}\rh{i}\rh{i+1} \si{i}= \si[2]{i} \cdot \sii{i}\rh{i}\sii{i}\cdot \si{i}\rh{i+1}\si{i}=\si[2]{i} \rh{i+1}\rh{i}$. Applying \repr{present}, we obtain:
\begin{align*}
\si{i}(\rh{n-1} \cdots \rh{1}) \sii{i} &= \rh{n-1}\cdots \rh{i+2} \si{i} \rh{i+1}\rh{i}\sii{i} \rh{i-1} \cdots \rh{1} = \si[2]{i} (\rh{n-1} \cdots \rh{1}).
\end{align*}
Hence:
\begin{equation}\label{eq:anti}
 \si{i}(\rh{n-1} \cdots \rh{1})= (\rh{n-1} \cdots \rh{1}) \sii{i}.
\end{equation}
and so $\si{i}$ commutes with $(\rh{n-1} \cdots \rh{1})^2$. Now let us show that  $\rh{n}$ commutes with $(\rh{n-1} \cdots \rh{1})^2$. Set $w=B_{1,n} \cdots B_{n-2,n} B_{n-1,n}$. Using \relem{commute} and \repr{present}, we see that
\begin{align}
(\rhi{1}\cdots \rhi{n-1})\rh{n} (\rh{n-1} \cdots \rh{1}) &= (\rhi{1}\cdots \rhi{n-2})\rh{n} (\rh{n-2} \cdots \rh{1})B_{n-1,n}\notag\\
&= (\rhi{1}\cdots \rhi{n-3})\rh{n} (\rh{n-3} \cdots \rh{1})B_{n-2,n} B_{n-1,n}\notag\\
&=\rh{n}w.\label{eq:rhw}
\end{align}
The relation $w= \si{n-1}\cdots \si{2}\si[2]{1} \si{2}\cdots \si{n-1}$ holds in $B_n(\rp)$ since it holds in $B_n(\dt)$. Thus $(\rhi{1}\cdots \rhi{n-1}) w (\rh{n-1} \cdots \rh{1})=w^{-1}$ using \req{anti}, and so from \req{rhw}, we obtain
\begin{align*}
(\rhi{1}\cdots \rhi{n-1})^2\rh{n} (\rh{n-1} \cdots \rh{1})^2&=(\rhi{1}\cdots \rhi{n-1})\rh{n}w (\rh{n-1} \cdots \rh{1})\\
&= \rh{n}ww^{-1}=\rh{n},
\end{align*}
which shows that $\rh{n}$ commutes with $(\rh{n-1} \cdots \rh{1})^2$.

Hence $(\rh{n-1} \cdots \rh{1})^2$ belongs to $Z(B_n(\rp))$, which is the cyclic subgroup of order~$2$ generated by $\ft$. Suppose that $(\rh{n-1} \cdots \rh{1})^2=1$. Clearly, $\rh{n-1} \cdots \rh{1}\neq 1$ (since under the projection $P_n(\rp)\to P_1(\rp)$, $\rh{n-1} \cdots \rh{1}$ is sent to $\rh{1}$). So $\rh{n-1} \cdots \rh{1}$ is an element of $P_n(\rp)$ of order~$2$. By \repr{fullpn}, it follows that $\rh{n-1} \cdots \rh{1}=\ft$. But this cannot be so, because $\ft$ is central in $B_n(\rp)$, while by \req{anti}, $\rh{n-1} \cdots \rh{1}$ is not. We conclude that $(\rh{n-1} \cdots \rh{1})^2=\ft$.
\end{proof}

Since by \repr{rhosq} the element $\rh{n-1}\cdots \rh{1}$ is of order~$4$, the characterisation of the torsion of $P_n(\rp)$ thus follows from \reco{tor24}:

\begin{cor}\label{cor:torspn}
Let $n\geq 2$. Then the (non-trivial) torsion of $P_n(\rp)$ is precisely~$2$ and~$4$.\hfill$\sq$
\end{cor}

We have a result analogous to that of \repr{rhosq} for $\rhs{n}{1}$:
\begin{prop}\label{prop:rhonsq}
Let $n\in\N$. In $B_n(\rp)$, $(\rhs{n}{1})^2=\ft$.
\end{prop}

\begin{proof}
If $n=1$ then $\rh[2]{1}=\ft[1]=1$. So suppose that $n\geq 2$. First, using \repr{present} and the definition of $w$ given in the proof of \repr{rhosq}, we see that
\begin{align}
\rh[2]{n} &= \sii{n-1}\cdots \sii{1} \rh{1} \sii{1}\cdots \sii{n-1} \cdot \sii{n-1}\cdots \sii{1} \rh{1} \sii{1}\cdots \sii{n-1}\notag\\
& = \sii{n-1}\cdots \sii[2]{1}\cdots \sii{n-1}= w^{-1}.\label{eq:rh2n}
\end{align}
By \repr{present} and \relem{commute}, we observe that:
\begin{align*}
(\rh{n}\cdots \rh{1})^ 2 &= \rh{n} \rh{n-1}\cdots \rh{1} \rh{n-1}\rh{n}\rh{n-2}\cdots \rh{1} B_{n-1,n}\\
&= \rh{n} (\rh{n-1}\cdots \rh{1})^2 \rh{n} w=\rh{n}\ft \rh{n} w=\ft \rh[2]{n}w=\ft,
\end{align*}
using \req{rh2n} and the fact that $\ft$ is central. The result follows.
\end{proof}

Now we go on to study $P_3(\rp)$ in more detail, using the fact that
it is a semi-direct product. Since $P_2(\rp)\cong Q_8$, $P_2(\rp)$
has~$6$ elements of order~$4$ which make up three distinct
conjugacy classes. We will show that the same happens in  $P_3(\rp)$.

\begin{prop}
Let $\map {p_{\ast}}{P_3(\rp)}{P_2(\rp)}$ be the homomorphism induced by the projection onto the first two coordinates. Then:
\begin{enumerate}
\item \label{it:surj4} $p_{\ast}$ sends the set of elements of
$P_3(\rp)$ of order~$4$ surjectively onto the set of elements of
$P_2(\rp)$ of order~$4$; 
\item \label{it:4con} any two elements of $P_3(\rp)$ of order~$4$ which project onto the same element  of $P_2(\rp)$ are conjugate; 
\item \label{it:4class} the set of elements of $P_3(\rp)$ of
order~$4$ form three distinct conjugacy classes;
\item \label{it:subiso} if $H$ is a subgroup of $P_3(\rp)$ abstractly
isomorphic to $Q_8$, then $p_{\ast}(H)\cong P_2(\rp)$. Further, any
two such subgroups are conjugate. 
\end{enumerate}
\end{prop}

\begin{proof}
The proof of this proposition will rely on the results of
\resec{rp2}, notably the semi-direct product $P_3(\rp)\cong \F\rtimes
P_2(\rp)$ given by \repr{p2p3}, and the corresponding representation
$\map{\phi}{P_2(\rp)}{\aut{\F}}$ of \reco{repres}. If $z\in
P_2(\rp)$, let us write $\phi(z)=\phi_z$. In what follows we shall
identify $P_3(\rp)$ with $\F\rtimes P_2(\rp)$, in particular, if
$(w,t)\in \F\rtimes P_2(\rp)$ then $p_{\ast}((w,t))=t\in P_2(\rp)$.

part~{(\ref{it:surj4})} follows directly from the existence of
the splitting. For part~{(\ref{it:4con})}, suppose that the 
element $(w,t)\in \F\rtimes P_2(\rp)$ is of order~$4$. By
part~{(\ref{it:surj4})}, $t\in P_2(\rp)$ is of order~$4$, and
hence so is $(1,t)\in \F\rtimes P_2(\rp)$. To prove the result, it
suffices to show that $(w,t)$ and $(1,t)$ are conjugate. From the
explicit description of $\phi$ given in \reco{repres}, one sees that
$\phi_{t^2}=\id_{\F}$, and so:
\begin{equation*}
(1,1)=(w,t)^4=(w\phi_t(w),t^2)^2 =
(w\phi_t(w)\phi_{t^2}(w\phi_t(w)),t^4) =((w\phi_t(w))^2,1).
\end{equation*}
Since $\F$ is torsion free and $(w\phi_t(w))^2=1$, it follows that 
$w\phi_t(w)=1$, and so $\phi_t(w)=w^{-1}$. Let us characterise such
$w\in \F$. Since $\phi_{t^2}=\id_{\F}$ for all $t\in P_2(\rp)$, we
see that $\phi_{t^3}=\phi_t$, and so it suffices to suppose that
$t\in\brak{\rh{1},\rh{2},\rh{1}\rh{2}}$. Using the description of the
associated automorphisms $\phi_t$ given in \reco{repres}, we see in
each of the three possible cases  that $w$ is of the (reduced) form
$w=w_1w_2$, where $w_2=\phi_t(w_1^{-1})$. This implies that:
\begin{equation*}
(w_1,1)\cdot (1,t)\cdot (w_1,1)^{-1}=(w_1,1)\cdot (1,t)\cdot (w_1^{-1},1)= (w_1,1)\cdot (\phi_t(w_1^{-1}),t)=(w,t),
\end{equation*}
and part~{(\ref{it:4con})} follows. One may then check easily
that if $(1,t)$ and $(1,t')$ are conjugate in $\F\rtimes P_2(\rp)$
then $t$ and $t'$ are conjugate in $P_2(\rp)$. Since the set of
elements of $P_2(\rp)$ of order~$4$ splits into three conjugacy
classes, part~{(\ref{it:4class})} follows from
parts~{(\ref{it:surj4})} and~{(\ref{it:4con})}.

For the first part of~{(\ref{it:subiso})}, since $H$ is a
torsion subgroup and the kernel of $p_{\ast}$ is torsion free, it
follows that $\map{p_{\ast}\vert_H}{H}{P_2(\rp)}$ is injective, and
thus an isomorphism.  For the second part, let $K$ be the subgroup of
$P_3(\rp)$ generated by $s(\rh{1})$ and $s(\rh{2})$, where $\map
{s}{P_2(\rp)}{P_3(\rp)}$ is the section of \repr{p2p3}, and let
$H$ be a subgroup of $P_3(\rp)$ which is isomorphic to $Q_8$. To
prove the statement, it suffices to prove that $H$ is conjugate to
$K$. From above, $p_{\ast}(H)=P_2(\rp)$, so let $b\in H$ be such that
$p_{\ast}(b)=\rh{1}$. By part~{(\ref{it:4con})}, there exists
$\gamma\in P_3(\rp)$ such that $\gamma b \gamma^{-1}=s(\rh{1})$. Let
$H'=\gamma H \gamma^{-1}$. Then it suffices to prove that $H'=K$. To
do so, first notice that $s(\rh{1})\in H'\cap K$. Let $c\in H'$ be such
that $p_{\ast}(c)=\rh{2}$. Again by part~{(\ref{it:4con})}, since
$p_{\ast}(s(\rh{2}))=\rh{2}$, there exists $(w,t)\in \F\rtimes P_2(\rp)$
such that:
\begin{equation}\label{eq:cconj}
c=(w,t)\cdot (1,\rh{2})\cdot (w,t)^{-1}=(w\phi_{t\rh{2}t^{-1}}(w^{-1}),t\rh{2}t^{-1}).
\end{equation}
We have once more identified $P_3(\rp)$ with $\F\rtimes P_2(\rp)$, in
particular, $s(\rh{i})=(1,\rh{i})$ for $i=1,2$. By \reco{repres}, we obtain
$\phi_{\rh{2}}=\phi_{\rhi{2}}$, and since $t\rh{2}t^{-1}$ is of the form
$\rh[\pm 1]{2}$,  it follows that $\phi_{t\rh{2}t^{-1}}= \phi_{\rh{2}}$. From
\req{cconj}, we obtain:
\begin{equation}\label{eq:cmone}
c^{-1}=(w\phi_{\rh{2}}(w^{-1}), \rh[\mp 1]{2}).
\end{equation}
On the other hand, $s(\rh{1})\cdot c \cdot (s(\rh{1}))^{-1}=c^{-1}$ in $H'$, so:
\begin{equation}\label{eq:cmtwo}
c^{-1} = (1,\rh{1})\cdot (w\phi_{\rh{2}}(w^{-1}),\rh[\pm 1]{2})\cdot (1,\rh{1})^{-1} = (\phi_{\rh{1}}(w\phi_{\rh{2}}(w^{-1})), \rh[\mp 1]{2}).
\end{equation}
Combining equations~\reqref{cmone} and~\reqref{cmtwo}, we see that
$w\phi_{\rh{2}}(w^{-1})$ is a fixed element of $\phi_{\rh{1}}$, but
by \reco{repres}, the only such element is the identity. So
$\phi_{\rh{2}}(w^{-1})=w^{-1}$, and similarly we conclude that $w=1$.
By \req{cconj}, we see that $c=(1,\rh[\pm 1]{2})=s(\rh[\pm 1]{2})$,
and since $H'$ is generated by $\brak{c,s(\rh{1})}$, it follows that
$H'=K$ as required.
\end{proof}

\begin{rem}\label{rem:mura}
For the case $P_3(\rp)$, our results are compatible with those of
Murasugi~\cite{M}. His classification shows that the torsion elements
of $B_3(\rp)$ are, up to powers and conjugacy, $\sigma_1\sigma_2$
(order~$6$), $\sigma_1\sigma_2\sigma_1$ and $\rh{1}\sigma_2$ (both of
order~$4$), $\rh{2}\sigma_1$ (order~$8$), and
$\rh{3}\sigma_2\sigma_1$ (order~$12$). The torsion elements of
$P_3(\rp)$ are, up to powers and conjugacy, $\ft[3]$ (order~$2$),
$\rh{2}\rh{1}$ and $\rh{3}\rh{2}\rh{1}$ (both of order~$4$). Then
$\sigma_2(\rh{2}\rh{1})^{\pm 1}\sigma_2^{-1}=
(\rhi{3}\rhi{2}\rh{3}\rh{2}\rh{3}\rh{1})^{\pm 1}$ projects onto
$\rh[\pm 1]{1}$, $(\rh{3}\rh{2}\rh{1})^{\pm 1}$ projects onto
$(\rh{2}\rh{1})^{\pm 1}$, and the element $\sigma_1\sigma_2(\rh{2}\rh{1})^{\pm 1}\sigma_2^{-1}\sigma_1^{-1}= (\rh[2]{1}\rhi{1}\rhi{2}\rh{1}\rh{2}\rh{3} \rhi{2}\rhi{1}\rh{2}\rh{1}\rh{2})^{\pm 1}$ projects onto $\rh[\mp
1]{2}$.
\end{rem}

\section{Roots of the full twist and torsion elements in $B_n(\rp)$}\label{sec:rootbn}

In this section, we study the torsion of the group $B_n(\rp)$ for
$n\geq 2$. In~\cite{M}, Murasugi introduces elements $A_1(n,r,s,q)$
and $A_2(n,r,s,q)$ of  $B_n(\rp)$, and shows that if $r>0$ then there
are the relations $mp=2rk$ and $mq=sk$ for integers $m\neq 0$, and
for these integers, we have $A_1(n,r,s,q)^{2rm}=\ft$ and
$A_2(n,r,s,q)^{4rm}=1$. But the order of such elements is not given.
We show that the full twist $\ft$ is the unique element of order two
in $B_n(\rp)$, and then prove \reth{torsbn} which says that
$B_n(\rp)$ has an element of order $k$ if and only if $k$ divides
either~$4n$ or~$4(n-1)$. As a consequence, we obtain \reth{rootbn},
so that the full twist $\ft$ has a $k\th$ root if and only if $k$
either divides~$2n$ or~$2(n-1)$. Our techniques are different to
those of Murasugi, in particular, we use the notion of intermediate
coverings between $F_n(\rp)$ and $D_n(\rp)$ introduced
in~\cite{GoG2}.

\begin{prop}\label{prop:unibn}
If $n\geq 2$ then the full twist $\ft$ is the unique element of  $B_n(\rp)$ of order~$2$.
\end{prop}

\begin{proof}
Let $\map{\pi}{B_n(\rp)}{\sn}$ denote the natural epimorphism, and
suppose that $\alpha$ is an element of $B_n(\rp)$ of order~$2$
different from $\ft$. So $\alpha\notin P_n(\rp)$ by \repr{fullpn}. It
follows that the cycle decomposition of $\pi(\alpha)$ is a product of
a non-zero number of transpositions and a certain number (perhaps
zero) of cycles of length one. Let $D_{2,n-2}(\rp)=
F_n(\rp)/\sn[2]\times \sn[n-2]$ be the intermediate covering between
$F_n(\rp)$ and $D_n(\rp)$~\cite{GoG2}, and let $B_{2,n-2}(\rp)$ denote
the corresponding subgroup of $B_n(\rp)$. Consider the homomorphism
$\map{\xi}{B_{2,n-2}(\rp)}{B_2(\rp)}$ induced by the projection
$D_{2,n-2}(\rp)\to D_2(\rp)$ onto the first two coordinates. By
conjugating $\alpha$ if necessary, we may suppose without loss of
generality that $\alpha\in B_{2,n-2}(\rp)$, and that the permutation
of $\xi(\alpha)$ is a transposition. Since $\alpha$ is of order~$2$
then so is $\xi(\alpha)$. But from~\cite{vB}, the group $B_2(\rp)$ is
known to be the dicyclic group of order~$16$ comprised of the
identity, $10$~elements of order~$4$, $4$~elements of order~$8$ and a
single element of order~$2$ which is the full twist $\ft[2]$. But
this cannot be the case because $\ft[2]\in P_2(\rp)$, so its
permutation is trivial. We thus obtain a contradiction, and this
proves the result.
\end{proof}    

We now derive a necessary condition on the permutation associated
with a torsion element.

\begin{lem}[Fundamental Lemma]\label{lem:fund}
Let $n\in\N$. Let $\alpha \in B_n(\rp)\setminus P_n(\rp)$ be a
torsion element. If the permutation $\pi(\alpha)$ contains at least
two cycles of length one then $\alpha\neq \ft$ but $\alpha^2=\ft$.
\end{lem}

\begin{proof}
Let $\alpha \in B_n(\rp)\setminus P_n(\rp)$ (so $n\geq 2$) be a
torsion element of order~$k$ such that the permutation $\pi(\alpha)$
contains at least two cycles of length one. Since $\alpha\notin
P_n(\rp)$, \repr{fullpn} implies that $\alpha\neq\ft$. Now $\ft$
belongs to $Z(B_n(\rp))$, so $\alpha$ would satisfy the
equation $\alpha^2=\ft$ if and only if its conjugates did too.
Conjugating $\alpha$ if necessary, we may thus assume without loss of
generality that $\alpha$ belongs to the subgroup $B_{1,1,n-2}(\rp)$
of $B_n(\rp)$ corresponding to the intermediate covering
$D_{1,1,n-2}(\rp)= F_n(\rp)/\sn[1]\times \sn[1]\times \sn[n-2]$.
Consider the homomorphism $\map{\xi}{B_{1,1,n-2}(\rp)}{P_2(\rp)}$
induced by the projection $D_{1,1,n-2}(\rp)\to F_2(\rp)$ onto the
first two coordinates. Now $P_2(\rp)\cong Q_8$, thus
$(\xi(\alpha))^2=\ft[2]^{\epsilon}$, where $\epsilon\in\brak{0,1}$.
Since $\ft\in B_{1,1,n-2}(\rp)$ and $\xi(\ft)=\ft[2]$, we have that
$\alpha^2\cdot \ft^{-\epsilon} \in\ker{\xi}=
B_{n-2}(\rp\setminus\brak{x_1,x_2})$ which we know to be torsion
free. But $(\alpha^2\cdot \ft^{-\epsilon})^{2k}=1$ because $\ft$ is
central in $B_n(\rp)$. Hence $\alpha^2=\ft^{\epsilon}$. If
$\epsilon=0$ then $\alpha^2=1$ which since $\alpha\neq 1$ implies
that $\alpha=\ft$ by \repr{unibn}. But we have already established
that this cannot be the case. Hence $\epsilon=1$, and thus
$\alpha^2=\ft$ as required.
\end{proof}

\begin{cor}\label{cor:torrp}
Let $n\in\N$. If $\alpha$ is a torsion element of $B_n(\rp)$ then its order divides $4n$ or $4(n-1)$. 
\end{cor}

\proof
Let $\ell$ denote the order of $\alpha$. If $\alpha\in P_n(\rp)$ then the result follows from \repr{fullpn}. So suppose that the cycle decomposition of $\pi(\alpha)$ contains at least one cycle of length greater than or equal to two. Then there are three possibilities:
\begin{enumerate}
\item the non-trivial cycles in the cycle decomposition of
$\pi(\alpha)$ are all of the same length $\ell_1$. Set $\ell_0=1$,
and for $i=0,1$, let $k_i$ be the number of cycles of length $\ell_i$
in the cycle decomposition of $\pi(\alpha)$. If $k_0\geq 2$ then it
follows from the Fundamental Lemma that $\alpha^2=\ft$, and hence
$\alpha$ is of order~$4$ which divides $4n$ as required. So suppose
that $k_0\in\brak{0,1}$. Then $\alpha^{\ell_1}\in P_n(\rp)$, and so
by \reco{tor24}, either $\alpha^{\ell_1}$ is the identity element,
and so $\alpha$ is of order $\ell$, or $\alpha^{\ell_1}$ is of
order~$2$, and so $\alpha$ is of order $2\ell$, or else
$\alpha^{\ell_1}$ is of order~$4$, and so $\alpha$ is of order
$4\ell$ (we have used the fact that $\alpha^j\notin P_n(\rp)$ for all
$1\leq j<\ell$). Since $n=k_0+\ell_1k_1$, we see that $\ell_1$
divides $n$ or $n-1$, and so the order of $\alpha$ divides $4n$ or
$4(n-1)$ as required.
\item the non-trivial cycles in the cycle decomposition of
$\pi(\alpha)$ are of two different lengths $\ell_1,\ell_2$ with
$1<\ell_1<\ell_2$. Set $\ell_0=1$, and for $i=0,1,2$, let  $k_i$ be
the number of cycles of length $\ell_i$ in the cycle decomposition of
$\pi(\alpha)$. Applying the Fundamental Lemma to $\alpha$, we see
that $k_0\in\brak{0,1}$. Consider the element $\alpha^{\ell_1}$;  it
belongs to $B_n(\rp)\setminus P_n(\rp)$, and contains at least two
cycles of length one. Applying the Fundamental Lemma, it follows that
$\ell_2=2\ell_1$, and $\alpha^{2\ell_1}=\ft$. Hence $\alpha$ is of
order $4\ell_1$. Now $n=k_0+k_1\ell_1+k_2\ell_2=
k_0+k_1\ell_1+2k_2\ell_1$, hence $\ell_1$ divides $n$ or $n-1$, and
thus the order $4\ell_1$ of $\alpha$ divides $4n$ or $4(n-1)$ as
required.
\item the cycle decomposition of $\pi(\alpha)$ contains cycles of at
least three different lengths $\ell_1,\ell_2,\ell_3$ with
$1<\ell_1<\ell_2<\ell_3$. Then $\alpha^{\ell_1}$ is a torsion element
belonging to $B_n(\rp)\setminus P_n(\rp)$ containing at least two
cycles of length one. By the Fundamental Lemma,
$\alpha^{2\ell_1}=\ft$, so the cycle decomposition of $\pi(\alpha)$
contains only transpositions and trivial cycles, but this is
impossible since $\ell_1<\ell_2<\ell_3$.\qed
\end{enumerate}

Let $n\geq 2$. We will show that there are elements in $B_n(\rp)$ of order $4n$ and of order $4(n-1)$. In $B_n(\rp)$, consider the following elements:
\begin{equation}\label{eq:abdef}
\left.
\begin{aligned}
a &= \sii{n-1}\cdots\sii{1}\cdot\rh{1}\\
b &= \sii{n-2}\cdots\sii{1}\cdot\rh{1}.
\end{aligned}
\right\}
\end{equation}

\begin{prop}\label{prop:a4n}
Let $n\geq 2$. In $B_n(\rp)$, $a$ is of order $4n$ and $b$ is of order $4(n-1)$.
\end{prop}

\begin{proof}
We start by proving that $a^{2n}=\ft$. As in \repr{rhosq}, we shall do this by showing that $a^{2n}$ 
is a non-trivial element of $B_n(\rp)$ that commutes with each of the generators
$\si{1},\ldots, \si{n-1},\rh{1}$.  We first establish two useful identities. 
Using \repr{present} and equations~\reqref{rhj1} and~\reqref{rh2n}, we have that:
\begin{align}\label{eq:rho1}
\sii{n-1}\cdots \sii{1} \rh{1}\si{1}\cdots \si{n-1}&= \sii{n-1}\cdots \sii{1}\rh{1} \sii{1}\cdots \sii{n-1} 
 \si{n-1}\cdots \si{2}\si[2]{1}\si{2} \cdots \si{n-1}\notag\\
&=\rh{n}\cdot \rhi[2]{n}=\rhi{n},
\end{align}
and
\begin{equation}\label{eq:rho12}
\si{1}\si{2} \cdots \si{n-1} \rh{1} \sii{n-1}\cdots \sii{1}= \si{1}\rh{1}\sii{1}=\si[2]{1} \rh{2}.
\end{equation}
Using \repr{present} and induction, we see that $a^{-i}\si{1} a^i= \si{i+1}$ for all $0\leq i\leq n-2$. Further, writing $a_1=\sii{n-1}\cdots \sii{1}$, we see that 
\begin{equation*}
a^{-1}\si{n-1} a =\rhi{1}\si{1}\cdots \si{n-2}\si{n-1}\si{n-1} \sii{n-1}\sii{n-2}\cdots \sii{i}\rh{1},
\end{equation*}
and applying the braid relations, it follows that $a^{-1}\si{n-1} a=\rhi{1}a_1\si{1} 
a_1^{-1}\rh{1}$. From this, we obtain $a^{-1}(a^{-1}\si{n-1} a)a=(a^{-1}\rhi{1}a_1a)\si{2} 
(a^{-1}\rhi{1}a_1a)^{-1}$. Equations~\reqref{rho1} and~\reqref{rho12} yield: 
\begin{align*}
a^{-1}\rhi{1}a_1a &= \rhi{1}\si{1}\cdots \si{n-1}\rhi{1}\sii{n-1}\cdots \sii{1}  \sii{n-1}\cdots \sii{1} \rh{1}\\&=\rhi{1}\rhi{2}\sii[2]{1} \rhi{n}\sii{n-1}\cdots \sii{1}.
\end{align*}
Since $\sii{n-1}\cdots \sii{1}\si{2}\si{1}\cdots \si{n-1}=\si{1}$, we  obtain
\begin{align*}
a^{-2}\si{n-1} a^2&= \rhi{1}\cdot \rhi{2}\sii[2]{1}\cdot \rhi{n}\sii{n-1}\cdots \sii{1}\cdot 
\si{2}\cdot \si{1}\cdots \si{n-1} \rh{n}\si[2]{1}\rh{2} \rh{1}\\
&=  \rhi{1}\rhi{2}\si{1} \rh{2} \rh{1}.
\end{align*}
Finally, by \req{anti}, we see that
\begin{equation*}
a^{-1}(a^{-1}\si{n-1} a)a=  \rhi{1}\rhi{2}\si{1} \rh{2} \rh{1}= \rhi{1}\cdots \rhi{n}\si{1}\rh{n}
\cdots \rh{1}=\sii{1}.
\end{equation*}
Thus the conjugates of $\si{1}$ by successive powers of $a^{-1}$ are:
\begin{equation*}
\si{1},\ldots, \si{n-1}, a^{-1}\si{n-1}a, \sii{1},\ldots, \sii{n-1}, a^{-1}\sii{n-1}a,
\end{equation*}
and then $\si{1}$ after $2n$ conjugations. In particular, $a^{-2n}\si{i} a^{2n}=\si{i}$ 
for all  $1\leq i\leq n-1$, and so $a^{2n}$ commutes with $\si{i}$. 

Now consider the successive conjugates of $\rh{1}$ by $a^{-1}$. Using \repr{present} and induction, we see that $a^{-i}\rh{1}a^i=\rh{i+1}$ for all $0\leq i\leq n-1$. By 
equations~\reqref{rhj1} and~\reqref{rh2n}, we have that
\begin{equation*}
\rh{n}= \sii{n-1}\cdots \sii{1} \rh{1} \sii{1}\cdots \sii{n-1}= a\rh{1} a^{-1}\rh[2]n.
\end{equation*}
It follows from this that $a^{-1}\rh{n} a=\rhi{1}$, and hence that the successive conjugates of 
$\rh{1}$ by $a^{-1}$ are $\rh{1},\ldots \rh{n}, \rhi{1},\ldots, \rhi{n}$ and then $\rh{1}$. So 
$a^{-2n}\rh{1}a^{2n}=\rh{1}$, and with the conclusion of the previous paragraph, we see that 
$a^{2n}$ is central in $B_n(\rp)$.

Suppose that $a^{2n}=1$. Since the permutation associated with $a$ is an $n$-cycle, it follows 
that $a$ is of order $n$ or $2n$, and so $a^n$ is central. If $a$ is of order $n$ this is 
obvious, while if $a$ is of order $2n$, we have that $(a^n)^2=1$, and so $a^n=\ft$ by
\repr{fullpn}. But $a^{-n}\si{1}a^n=\sii{1}$ from above, which yields a contradiction. 
Since $Z(B_n(\rp))$ is generated by the full twist $\ft$ of order~$2$, we deduce that 
$a^{2n}=\ft$, and since $n$ divides the order of $a$, it follows that $a$ is of order $4n$.

Now let us show that $b^{2(n-1)}=\ft$. Notice first that $b=\si{n-1}a$. By induction, we see that for 
all $1\leq j\leq n-1$, $b^j=a^j \sii{j-1}\cdots \sii{1}\cdot a^{-1}\si{n-1} a$. Indeed if the result is true for $1\leq j\leq n-2$ then using the conclusion of the first paragraph of this proof, we see that
\begin{align*}
b^{j+1} &= a^j \sii{j-1}\cdots \sii{1}\cdot a^{-1}\si{n-1} a\cdot \si{n-1}a\\
&=
a^{j+1} a^{-1} \sii{j-1}\cdots \sii{1}\cdot a a^{-2}\si{n-1} a^2\cdot a^{-1}\si{n-1}a\\
&= 
a^{j+1} \sii{j}\cdots \sii{1} a^{-1}\si{n-1}a.
\end{align*}
Hence
\begin{align*}
b^{2(n-1)}  &= (a^{n-1} \sii{n-2}\cdots \sii{1} a^{-1}\si{n-1}a)^2\\
&= a^{2n-(n+1)}(\sii{n-2}\cdots \sii{1} a^{-1}\si{n-1}a)a^{n+1-2} 
(\sii{n-2}\cdots \sii{1} a^{-1}\si{n-1}a)a a^{-1}\\
& = \ft \si{n-1}\cdots \si{1} a^{-1} 
\sii{n-1}\cdots \sii{1}a^{-1}\\
&= \ft \si{n-1}\cdots \si{1} \rhi[2]{1} \si{1}\cdots \si{n-1}= \ft,
\end{align*}
by \repr{present}. As for the case of $a$, to see that $b$ is of order~$4(n-1)$, it suffices to 
prove that  $b^{n-1}$ is not central. An easy induction argument shows that 
$b^{-i}\si{1} b^i=\si{i+1}$ for $0\leq i\leq n-3$. Further, a calculation analogous to that above of $a^{-2}\si{n-1}a^2$ shows that $b^{-2}\si{n-1}b^2=\sii{1}$. So $b^{-(n-1)} \si{1} b^{(n-1)}=\sii{1}$, thus $b^{n-1}$ is not central, and the proof of the proposition is complete.
\end{proof}

This enables us to prove \reth{torsbn}:

\begin{proof}[Proof of \reth{torsbn}]
To prove the `only if' part, it suffices to consider an appropriate
power of one of the elements $a$ and $b$, and to apply \repr{a4n}.
The `if' part follows directly from \reco{torrp}.
\end{proof}

\begin{rem}
It is possible (by an induction argument for example) to prove that $b^{n-1}=\rh{n-1}\cdots \rh{1}$ and $a^n=\rh{n}\cdots \rh{1}$. Then the `only if' part of \reth{torsbn} also follows from Propositions~\ref{prop:rhosq} and~\ref{prop:rhonsq}.
\end{rem}

Finally, we are able to prove \reth{rootbn}.

\begin{proof}[Proof of \reth{rootbn}]
To prove the `only if' part, let $\alpha$ be a $k\th$ root of the
full twist. Then $\alpha$ is a torsion element of $B_n(\rp)$ of order
$2k$. By \reth{torsbn}, it follows that $2k$ divides either $4n$ or
$4(n-1)$, and the result follows. To prove the `if part', the element
$a$ (respectively, $b$) given in \req{abdef} has the property that
its $2n\th$ (respectively, $2(n-1)\th$) power is of order~$2$, and so
is the full twist $\ft$ by \repr{unibn}. The result then follows by
taking an appropriate power of $a$ or $b$.
\end{proof}  

\begin{rem}
In the case $n=3$, it follows from \rerem{mura} that our results in this section are compatible with those of Murasugi.
\end{rem}

%%%%%%%%%%%%%%%%%%%%   End of main body of article
%
%                             References
%

\Addresses\recd

\end{document}

%% file: agtout.tex
\def\ifplaintex{\expandafter\ifx\csname documentclass\endcsname\relax}

\def\gtp{{\mathsurround=0pt\it $\cal G\mskip-2mu$eometry \&\ 
$\cal T\!\!$opology $\cal P\!$ublications}}  % GT publications

\def\recd{{\small Received:\qua\receiveddate\ifx\reviseddate\relax
\else\qquad Revised:\qua\reviseddate\fi\par}} 

\def\lognumber#1{\def\thelognumber{#1}}
\def\volumenumber#1{\def\thevolumenumber{#1}}
\def\volumeyear#1{\def\thevolumeyear{#1}}
\def\papernumber#1{\def\thepapernumber{#1}}
\def\pagenumbers#1#2{\def\startpage{#1}\def\finishpage{#2}}
\def\published#1{\def\publishdate{#1}}

\def\received#1{\def\receiveddate{#1}}

\def\accepted#1{\def\accepteddate{#1}}

\def\asciiauthors#1{\def\theasciiauthors{#1}}
\def\asciiaddress#1{\def\theasciiaddress{#1}}
\def\asciiemail#1{\def\theasciiemail{#1}}

\def\coverauthors#1{\def\thecoverauthors{#1}}
\long\def\asciiabstract#1{\long\def\theasciiabstract{#1}}

\let\\\par\let\thelognumber\relax\let\thevolumenumber\relax
\let\thepapernumber\relax\let\thevolumeyear\relax\let\startpage\relax
\let\finishpage\relax\let\publishdate\relax\let\receiveddate\relax
\let\reviseddate\relax\let\accepteddate\relax\let\theasciititle\relax
\let\theasciiauthors\relax\let\theasciiaddress\relax
\let\theasciiabstract\relax

\let\thecoverauthors\relax\let\theasciiemail\relax

\ifplaintex
\font\logobig=cmssbx10 scaled 3836
\font\logomed=cmssbx10 scaled 2557
\else
\font\logobig=cmssbx10 scaled 4200
\font\logomed=cmssbx10 scaled 2800
\fi

\long\def\makeagttitle{   %%% start of definition of \makeagttitle
\count0=\startpage
\agt\hfill      %   Journal title (top left) 
\hbox to 45truept{\vbox to 0pt{\vglue -13truept{\logomed A\kern -.37em{\logobig 
T}\kern -.38em G}\vss}\hss}
\break
{\small Volume \thevolumenumber\ (\thevolumeyear)
\startpage--\finishpage\nl
Published: \publishdate}

\vglue .25truein

{\parskip=0pt\leftskip 0pt plus
1fil\def\\{\par\smallskip}{\Large\bf\thetitle}\par\medskip} \vglue
0.05truein

{\parskip=0pt\leftskip 0pt plus 1fil\def\\{\par}{\sc\theauthors}
\par\medskip}%
 
\vglue 0.03truein 

{\small\leftskip 25truept\rightskip 25truept{\bf Abstract}\stdspace\theabstract

{\bf AMS Classification}\stdspace\theprimaryclass
\ifx\thesecondaryclass\relax\else; \thesecondaryclass\fi\par
{\bf Keywords}\stdspace \thekeywords\par}\vglue 7truept

}   %%%% end of definition of \makeagttitle

\ifplaintex
\font\phead=cmsl9 scaled 950
\font\pnum=cmbx10 scaled 913
\font\pfoot=cmsl9 scaled 950
\headline{\vbox to 0pt{\vskip -4.5mm\line{\small\phead\ifnum
\count0=\startpage ISSN 1472-2739 (on-line) 1472-2747 (printed)
\hfill {\pnum\folio}\else\ifodd\count0\def\\{ }% 
\ifx\theshorttitle\relax\thetitle\else\theshorttitle\fi\hfill{\pnum\folio}
\else\def\\{ and }{\pnum\folio}\hfill\ifx\theshortauthors\relax\theauthors
\else\theshortauthors\fi\fi\fi}\vss}}
\footline{\vbox to 0pt{\vglue 0mm\line{\small\pfoot\ifnum\count0=\startpage
\copyright\ \gtp\hfill\else
\agt, Volume \thevolumenumber\ (\thevolumeyear)\hfill\fi}\vss}}
\else
\font=cmsl9 scaled 1050
\font\lnum=cmbx10 
\font=cmsl9 scaled 1050
\makeatletter
\def\@oddhead{{\small\lhead\ifnum\count0=\startpage ISSN 1472-2739 
(on-line) 1472-2747 (printed)\hfill {\lnum\number\count0}\else\ifodd\count0
\def\\{ }\ifx\theshorttitle\relax \thetitle \else\theshorttitle\fi\hfill
{\lnum\number\count0}\else\def\\{ and }{\lnum\number\count0}
\hfill\ifx\theshortauthors\relax 
\theauthors\else\theshortauthors\fi\fi\fi}}\def\@evenhead{\@oddhead}
\def\@oddfoot{\small\lfoot\ifnum\count0=\startpage\copyright\ \gtp\hfill\else
\agt, Volume \thevolumenumber\ (\thevolumeyear)\hfill\fi}
\def\@evenfoot{\@oddfoot}
\makeatother
\fi
\let\maketitlepage\makeagttitle
\let\makeshorttitle\maketitlepage
\let\maketitle\maketitlepage

\newwrite\gtoutfile
\long\gdef\makeheadfile{  %%% start of definition of \makeheadfile
{\def\\{, }\def\s{ }
\immediate\openout\gtoutfile head.xxx
\immediate\write\gtoutfile{Proxy-for: \ifx\theasciiauthors\relax
\theauthors\else\theasciiauthors\fi\s<\ifx\theasciiemail\relax\theemail\else\theasciiemail\fi>}
\immediate\write\gtoutfile{\noexpand\\}
\immediate\write\gtoutfile{Authors: \ifx\theasciiauthors\relax
\theauthors\else\theasciiauthors\fi}
{\def\\{ }\immediate\write\gtoutfile{Title: \ifx\theasciititle\relax
\thetitle\else\theasciititle\fi}}
\immediate\write\gtoutfile{Subj-class: GT or SG, GR etc}
\immediate\write\gtoutfile{MSC-class: \theprimaryclass\ifx\thesecondaryclass\relax\else, \thesecondaryclass\fi}
\immediate\write\gtoutfile{Journal-ref: Algebr. Geom. Topol. \thevolumenumber\s
(\thevolumeyear) \startpage-\finishpage}
\immediate\write\gtoutfile{Comments: Published by Algebraic and
Geometric Topology at}
\immediate\write\gtoutfile{\s\s\s  http://www.maths.warwick.ac.uk/agt/AGTVol\thevolumenumber/agt-\thevolumenumber-\thepapernumber.abs.html}
\immediate\write\gtoutfile{\noexpand\\}
\immediate\write\gtoutfile{}
\ifx\theasciiabstract\relax
\immediate\write\gtoutfile{\theabstract}\else
\immediate\write\gtoutfile{\theasciiabstract}\fi
\immediate\write\gtoutfile{}
\immediate\write\gtoutfile{\noexpand\\}
\immediate\write\gtoutfile{}
\immediate\closeout\gtoutfile}}  %%% end of definition of \makeheadfile

\def\maketitlepage{\makeagttitle\makeheadfile}
\let\makeshorttitle\maketitlepage
\let\maketitle\maketitlepage